\documentclass[letterpaper, 12pt]{article}
\usepackage{amssymb,amsmath}
\usepackage{epsf}
\usepackage{times,bm,mathrsfs}
\usepackage{amsmath}
\usepackage{amssymb}
\usepackage{amsfonts}
\usepackage{mathrsfs}
\usepackage{bbm}
\usepackage{color}
\usepackage{colordvi}
\usepackage{graphicx}
\usepackage{amscd}
\usepackage{epsfig}

\newtheorem{theorem}{Theorem}

\newtheorem{definition}{Definition}

\newtheorem{lemma}{Lemma}

\newtheorem{remark}{Remark}

\newtheorem{assumption}{Assumption}
\newenvironment{proof}{\noindent{\textbf{Proof:}}}{$\blacksquare$\vskip\belowdisplayskip}

\newenvironment{prevproof}[2]{\noindent {\bf {Proof ({#1}~\ref{#2}):}}}{$\blacksquare$\vskip \belowdisplayskip}

\definecolor{Red}{rgb}{1,0,0}
\definecolor{Blue}{rgb}{0,0,1}
\definecolor{Olive}{rgb}{0.41,0.55,0.13}
\definecolor{Green}{rgb}{0,1,0}
\definecolor{MGreen}{rgb}{0,0.8,0}
\definecolor{DGreen}{rgb}{0,0.55,0}
\definecolor{Yellow}{rgb}{1,1,0}
\definecolor{Cyan}{rgb}{0,1,1}
\definecolor{Magenta}{rgb}{1,0,1}
\definecolor{Orange}{rgb}{1,.5,0}
\definecolor{Violet}{rgb}{.5,0,.5}
\definecolor{Purple}{rgb}{.75,0,.25}
\definecolor{Brown}{rgb}{.75,.5,.25}
\definecolor{Grey}{rgb}{.5,.5,.5}
\definecolor{Black}{rgb}{0,0,0}

\newcommand{\ecal}{\mathcal{E}}

\newcommand{\gcal}{\mathcal{G}}

\newcommand{\pcal}{\mathcal{P}}
\newcommand{\qcal}{\mathcal{Q}}

\newcommand{\tcal}{\mathcal{T}}

\newcommand{\vcal}{\mathcal{V}}

\newcommand{\real}{\mathbb{R}}

\newcommand{\eps}{\varepsilon}

\newcommand{\prob}{\mathbb{P}}
\newcommand{\expec}{\mathbb{E}}
\newcommand{\E}{\mathbb{E}}
\renewcommand{\P}{\mathbb{P}}
\newcommand{\var}{\mathrm{Var}}

\newcommand{\median}{\mathrm{Median}}

\newcommand{\poly}{\mathrm{poly}}

\renewcommand{\root}{r}

\newcommand{\path}{\mathrm{P}}

\newcommand{\ball}[2]{\mathcal{C}_{#1}^{(#2)}}

\newcommand{\locations}{\mathcal{X}}
\newcommand{\mrca}{\mathrm{MRCA}}

\renewcommand{\time}{\tau}
\newcommand{\lgt}{\lambda}

\newcommand{\weight}{\omega}

\newcommand{\dist}{\mathrm{D}}

\newcommand{\distm}{\dist_\mathrm{m}}

\newcommand{\maxlgt}{\overline{\lgt}}
\newcommand{\minlgt}{\underline{\lgt}}

\newcommand{\lgtweight}{\Lambda}
\newcommand{\lgttotal}{\boldsymbol{\Lambda}}
\newcommand{\lgttotalextinct}{\boldsymbol{\Lambda}_{\mathrm{tot}}}
\newcommand{\mintime}{\underline{\time}}
\newcommand{\maxtime}{\overline{\time}}

\newcommand{\mrqw}{\boldsymbol{\Upsilon}^{(4)}}
\newcommand{\mrpw}{\boldsymbol{\Upsilon}^{(2)}}
\newcommand{\hw}{\mathbf{H}}
\newcommand{\hwfrac}{\gamma}
\newcommand{\hwgenes}{\mathbf{G}^H}
\newcommand{\maxhwfrac}{\overline{\hwfrac}}
\newcommand{\minhwfrac}{\underline{\hwfrac}}
\newcommand{\ratiolgt}{\rho_\lgt}
\newcommand{\ratiotime}{\rho_\time}

\newcommand{\remove}[1]{{}}

\begin{document}

\title{\vspace{0cm}
Recovering the tree-like trend
of evolution despite
extensive lateral genetic transfer:
A probabilistic analysis\footnote{
Preliminary results were announced without
proof in the proceedings of RECOMB 2012.
}
}

\author{
Sebastien Roch\footnote{Department of Mathematics
and Bioinformatics Program, UCLA.
Supported by NSF grant DMS-1007144.
This work
was done while SR was visiting the Institute for
Pure and Applied Mathematics (IPAM).}
\and
Sagi Snir\footnote{Institute of Evolution, 
University of Haifa. 
Supported by the USA-Israel Binational Science Foundation and by the Israel Science Foundation.
This work
was done while SS was visiting the Institute for
Pure and Applied Mathematics (IPAM).
}
}
\maketitle

\begin{abstract}
Lateral gene transfer (LGT) is a common mechanism of non-vertical
evolution where genetic material is transferred between two more or
less distantly related organisms.  It is particularly common in
bacteria where it contributes to adaptive evolution with important
medical implications.  In evolutionary studies, LGT has been shown to
create widespread discordance between gene trees as genomes become
mosaics of gene histories.  In particular, the Tree of Life has been
questioned as an appropriate representation of bacterial evolutionary
history.  Nevertheless a common hypothesis is that prokaryotic
evolution is primarily tree-like, but that the underlying trend is
obscured by LGT.  Extensive empirical work has sought to extract a
common tree-like signal from conflicting gene trees.
Here we give a probabilistic perspective on the problem of recovering
the tree-like trend despite LGT.  Under a model of randomly
distributed LGT, we show that the species phylogeny can be
reconstructed even in the presence of surprisingly many (almost
linear number of)  LGT events per gene tree. Our results, which are
optimal up to logarithmic factors, are based on the analysis of a
robust, computationally efficient reconstruction method and provides
insight into the design of such methods.  Finally we show that our
results have implications for the discovery of highways of gene
sharing.
\end{abstract}



\section{Introduction}\label{section:introduction}

High-throughput sequencing is transforming the study of evolution by allowing the integration of genome analysis and systematic studies, an area called phylogenomics~\cite{EisenFraser:03,DeBrPh:05}. An important step
in most
phylogenomic analyses is 
the reconstruction of a
tree of ancestor-descendant 
relationships---a gene tree---for 
each family of orthologous genes in a dataset. 
Such analyses have revealed widespread discordance between gene trees~\cite{GaltierDaubin:08},
leading some to question the meaningfulness
of the Tree of Life~\cite{GoDoLa:02,ZhLaGo:04,GogartenTownsend:05,BaSuLe+:05,DoolittleBapteste:07,Koonin:07}. In addition to statistical errors in gene tree estimation, various mechanisms commonly lead to incongruences between inferred gene histories, including hybridization events, 
duplications and losses in gene families, 
incomplete lineage sorting,
and lateral genetic transfers~\cite{Maddison:97}. 

Here we study specifically lateral
gene transfer (LGT), that is, the 
non-vertical transfer of genes
between more or less distantly related organisms
(as opposed to the standard vertical transmission
between parent and offspring). 
Estimates of the fraction of genes that have undergone
LGT vary widely---with some as high
as 99\%. See e.g.~\cite{DaganMartin:06,GaltierDaubin:08}
and references therein. LGT is particularly
common in bacterial evolution
and it has been recognized to 
play an important role in
microbial adaptation, selection and evolution
with implications in the study of infectious diseases~\cite{SmetsBarkay:05}.
As a result, the bacterial phylogeny 
is usually inferred from
genes that are thought to be immune to
LGT, typically ribosomal RNA genes.
However there is growing evidence 
that even such genes have in fact experienced 
LGT~\cite{YaZhWa:99,BeTePa+:03,ScScJa:03,DeShSc+:05}.
In any case, LGT appears to be a major source of
conflict between gene trees that must be taken into account
appropriately in phylogenomic analyses, in particular
when building phylogenies. 
This is the problem we address in this paper.

Despite the confounding
effect of LGT, we operate under the prevailing assumption 
that the evolution of organisms is
governed primarily by vertical inheritance. 
In particular we ask: 
\begin{enumerate}
\item 
How much genetic transfer
can be handled before 
the 
tree-like signal is completely erased? 
\item What phylogenetic reconstruction methods
are most effective under this
hypothesis? 
\end{enumerate}
These questions, and other related issues, 
have been the subject
of some empirical and simulation-based 
work~\cite{BeHaRa:05,GeWaKi:05,Galtier:07,PuWoKo:09,PuWoKo:10,KoPuWo:11}.
See also~\cite{GaltierDaubin:08,RaganBeiko:09}
for enlightening discussions.
In particular there is ample evidence that 
a strong tree-like signal can be extracted
in the presence of extensive 
LGT (although some debate
remains on this question~\cite{GoDoLa:02}).

In this paper we provide 
the first (to our knowledge) mathematical 
analysis of the issues above. 
We work under a 
stochastic 
model of gene tree topologies
positing that LGT events occur at more or less 
random locations
on the species phylogeny~\cite{Galtier:07}. 
In our main 
result we establish quantitative bounds implying 
that surprisingly
high levels of LGT---almost linear in the number of
branches for each gene---can be handled by
simple, computationally efficient inference procedures.  
That amount of genetic transfer appears to be
much higher than known empirical 
estimates of LGT frequency based on genomic datasets
in prokaryotes\footnote{Note that such estimates are typically based
on small numbers of genomes and, therefore,
are probably lower than reality~\cite{GaltierDaubin:08}.}.
Hence our results indicate that 
an accurate, reliable bacterial phylogeny 
should be reconstructible
if the vertical inheritance hypothesis is correct.
We prove that our bound on the achievable
rate of LGT is tight
up to logarithmic factors.
We also show that constraining LGT
to closely
related species makes the tree reconstruction
problem significantly easier.

Our theoretical approach complements simulation-based
studies in allowing a broad range of parameters
and tree shapes
to be considered.
Moreover our analysis provides new insights into the 
design of effective reconstruction methods
in the presence of
LGT. 
More precisely we focus on 
methodologies---both distance-based~\cite{KimSalisbury:01} 
and quartet-based~\cite{ZhGoCh+:06}---that
derive their statistical power from the aggregation of
basic topological information across genes.

In addition, we study the effect of
so-called highways of gene sharing,
roughly, preferred genetic exchanges
between {\em specific} groups of species.
Beiko et al.~\cite{BeHaRa:05}
provided empirical evidence
for the existence of such highways.
To identify highways,
they inferred LGT events 
by reconciling gene trees
with a trusted species tree.
In subsequent work, 
Bansal et al.~\cite{BaBaGo+:11}
formalized the problem
and designed a
fast highway detection algorithm 
that aggregates
conflicting signal across genes rather
than solving the difficult LGT inference problem
on each gene tree.
Similarly to Beiko et al., 
Bansal et al.~rely on a trusted species tree.

Here we show that a species
phylogeny can be reliably estimated in the 
presence of {\em both} random LGT events
and highways of LGT as long as such
highways involve a small enough fraction
of genes. Under extra assumptions, we also
design an algorithm for inferring
the location of highways.
Because we first recover the species phylogeny,
our highway reconstruction algorithm
does not require a trusted species tree.
In essence, our results on highways 
indicate that robust phylogeny
reconstruction in the presence of random LGT 
extends to a phylogenetic
network setting. For background on
phylogenetic networks, see e.g.~\cite{HuRuSc:10}.

We note that there exist related lines of work
in phylogenomics
addressing the issue of 
incomplete lineage sorting~\cite{DegnanRosenberg:09} 
in the presence of gene transfers and
hybridization events~\cite{ThRuIn+:07,JoMcLo:09,Kubatko:09,MengKubatko:09,YuThDe+:11,ChungAne:11}
as well as work on probabilistic models
involving gene duplications and losses~\cite{ArLaSe:09,CsurosMiklos:06}.

The rest of the paper is
organized as follows.
In Section~\ref{sec:results},
we define a stochastic model of LGT
and state our main results. 
A high-level description of our analysis
is given in Section~\ref{sec:analysis}.
Finally in Section~\ref{sec:highways}
we extend our results to
highways of gene sharing.

The results presented here were
announced without proof in~\cite{RochSnir:12}.

\section{Model and Main Results}\label{sec:results}

Before stating our main results, 
we present a stochastic model of LGT.
Roughly, following Galtier~\cite{Galtier:07},
we assume that LGT events occur more
or less at random along the species
phylogeny. Such a model appears to be consistent
with empirical evidence~\cite{GaltierDaubin:08}.

\paragraph{Notation} Recall that, for functions $f(n), g(n)$,
$f = O(g)$ means that there is constant $C > 0$ such that
$f(n) \leq C g(n)$ for all $n$ large enough. Similarly,
$f = \Omega(g)$ indicates $f(n) \geq C' g(n)$ for $C' > 0$. 
In addition $f = \Theta(g)$ is equivalent to
$f = O(g)$ and $f = \Omega(g)$.
By {\it polynomial in $n$}, we mean $O(n^{C''})$
for some constant $C'' > 0$.
We use the notation $\P[\ecal_0\,|\,\ecal_1]$
for the conditional probability of $\ecal_0$
given $\ecal_1$.

\subsection{Stochastic Model of LGT}

\paragraph{Gene trees and species phylogeny}
A {\em species phylogeny} (or
phylogeny for short) is a graphical representation 
of the speciation history of a group of organisms.
The leaves correspond to extant or extinct species. 
Each branching indicates a speciation event. 
Moreover we associate to each edge a 
positive value corresponding to the time
elapsed along that edge. 
For a tree $\tcal = (\vcal,\ecal)$ with leaf set $L$
and a subset of leaves $X \subseteq L$,
we let $\tcal|X$ be the {\em restriction of $\tcal$
to $X$}, that is, the subtree of $\tcal$ 
where we keep only those vertices and edges
on paths connecting two leaves in $X$. We say that $\tcal$ {\em
  agrees} (or is {\em consistent}) with $\tcal|X$.

\begin{definition}[Phylogeny]\label{def:phylogeny}
A \emph{(species) phylogeny} $T_s = (V_s,E_s,L_s;\root,\time)$ 
is a rooted tree with vertex set $V_s$, edge set $E_s$ and
$n$ (labelled) leaves $L_s = [n] = \{1,\ldots,n\}$ such that 1) the degree of all internal vertices
$V_s-L_s$ is exactly $3$ except the root $\root$
which has degree $2$, and 2) the edges are assigned 
inter-speciation times 
$\time : E_s \to (0,+\infty)$.
We assume that $T_s$ includes $n^+ > 0$ extant species
$L_s^+$
and $n^- \geq 0$ extinct species $L_s^-$, where
$n = n^+ + n^-$.
We also associate to each edge $e\in E_s$ in $T_s$ 
a {\em rate of lateral gene transfer} $0 < \lgt(e) < +\infty$. 
We denote by $T_s^+ = (V_s^+,E_s^+,L_s^+;\root,\time^+)$, 
the subtree of $T_s$
restricted to the extant leaves $L_s^+$,
that is, $T_s^+ = T_s|L_s^+$ rooted 
at the most recent common ancestor of $L_s^+$. 
We further suppress vertices of degree $2$ in $T_s^+$
except the root (in which case we add up the branch
lengths to obtain $\time^+$).
We call $T_s^+$ the {\em extant phylogeny}.
We assume that $T_s^+$ is ultrametric,
that is, from every node, 
the path lengths from that node
to all its descendant leaves are equal.
\end{definition}

Although we are ultimately interested in 
recovering the extant phylogeny, we include
extinct species in the model as they can be involved
in LGT events that affect the extant restriction
of the tree. See e.g.~\cite{Maddison:97}.

To infer the species phylogeny, 
we first reconstruct gene trees,
that is, trees of ancestor-descendant relationships 
for orthologous genes or loci.
Phylogenomic studies have revealed extensive 
discordance between such gene trees (e.g.~\cite{BaSuLe+:05,DoolittleBapteste:07}).
\begin{definition}[Gene tree]\label{def:genetree}
A \emph{gene tree} $T_g = (V_g,E_g,L_g;\weight_g)$
for gene $g$ 
is an unrooted tree 
with vertex set $V_g$, edge set $E_g$ and
$0 < n_g \leq n$ (labelled) leaves $L_g \subseteq \{1,\ldots,n\}$ with $|L_g| = n_g$ 
such that 1) the degree of every internal vertex
is either $2$ or $3$, and 2) the edges are assigned branch lengths $\weight_g : E_g \to (0,+\infty)$.
We let $\tcal_g = \tcal[T_g]$ be the
{\em topology} of $T_g$ where each internal
vertex of degree $2$ is suppressed.
\end{definition}
\begin{remark}[Gene trees vs. species phylogeny]
As we will discuss below, gene trees are derived from---
or ``evolve'' on---the species phlyogeny. They
may differ from the species phylogeny 
for various reasons. First, in our model,
their branch lengths represent expected
numbers of substitutions, instead time elapsed. 
Moreover, their topology may differ as a result, in our case, of
LGT events. See more details below.
\end{remark}
\begin{remark}[Rooted vs. unrooted]
Our stochastic model
of LGT requires a {\em rooted} species
phylogeny as time plays an important role
in constraining valid LGT events.
See, e.g.,~\cite{JIN-TCBB-2009}. 
In particular our
results rely on the ultrametricity property of the
extant phylogeny.
In contrast, branch lengths in gene trees correspond
to expected numbers of substitutions. 
As a result, gene trees are typically {\em unrooted}
and do not satisfy ultrametricity.
\end{remark}
\begin{remark}[Taxon sampling]
Each leaf in a gene tree corresponds to 
an extant species in the species phylogeny.
However, because of 
gene loss and taxon sampling, a taxon 
may not be represented in every gene tree.
\end{remark}
\begin{remark}[Branch lengths]
Each branch $e$ in a gene tree $T_g$
corresponds to a full or partial edge
in the species phylogeny $T_s$. 
In particular, we allow internal vertices
of degree $2$ in a gene tree to potentially
delineate between two consecutive species edges.
We allow the branch lengths $\weight_g(e)$
to be arbitrary, but one could easily
consider cases
where the branch lengths are
determined by inter-speciation times,
lineage-specific rates of substitution and
gene-specific rates of substitution.
The branch lengths will play 
a role in Section~\ref{sec:sequence}.
\end{remark}

\paragraph{Random LGT}
We formalize a stochastic model of LGT
similar to Galtier's~\cite{Galtier:07}. 
See also~\cite{KimSalisbury:01,Suchard:05,Jin-Bioinf-2006}
for related models.
The model accounts for LGT events
originating at random locations on the species
phylogeny with LGT rate 
$\lgt(e)$ prevailing along edge $e$. 

We will need the following notation.
Let $T_s = (V_s,E_s,L_s;\root,\time)$
be a fixed species phylogeny.
By a {\em location} in $T_s$, we mean
any position along $T_s$ seen as a continuous
object (also called $\real$-tree), that is,
a point $x$ along an edge $e \in E_s$.
We write $x \in e$ in that case.
We denote the set of locations in $T_s$
by $\locations_s$.
For any two locations $x, y$ in $\locations_s$,
we let $\mrca(x,y)$ be their most 
recent common ancestor (MRCA) in $T_s$
and we let $\time(x,y)$ be the length 
of the path connecting $x$ and $y$
in $T_s$ under the metric naturally defined
by the weights $\{\time(e), e \in E_s\}$,
interpolated linearly to locations
along an edge.
In words $\time(x,y)$, which we refer to
as the $\time$-distance between $x$ and $y$, 
is the sum of times
to $x$ and $y$ from $\mrca(x,y)$. We say that
two locations $x, y$ are 
{\em contemporaneous} if their respective
$\time$-distance to the root $\root$
is identical, that is, 
$$
\time(x,\root) = \time(y,\root).
$$ 
For $R > 0$,  we let 
$$
\ball{x}{R}
= \{y \in \locations_s : \time(\root,x) = \time(\root,y),
\ \time(x,y) \leq 2 R\}$$ 
be the
set of locations contemporaneous to $x$
at $\time$-distance at most $2R$ from $x$
(or in other words with MRCA at 
$\time$-distance at most $R$).
In particular, $\ball{x}{\infty}$ denotes the
set of all locations contemporaneous to $x$.
We let 
$\lgtweight(e) =  \lgt(e) \time(e)$, $e \in E_s$.
We note that, since $\lgt(e)$ is the LGT rate on $e$,
$\lgtweight(e)$ gives the expected number of 
LGT events along $e$. Further, we let
$$
\lgttotalextinct = \sum_{e\in E_s} \lgtweight(e),
$$
be the {\em total LGT weight} of the phylogeny 
and
$$
\lgttotal = \sum_{e\in \ecal(T_s|L_s^+)} \lgtweight(e),
$$
be the total LGT weight of the extant phylogeny,
where $\ecal(T_s|L_s^+)$ denotes the edge set
of $T_s|L_s^+$.

\begin{figure}[t]
\centering
\includegraphics[width=4in]{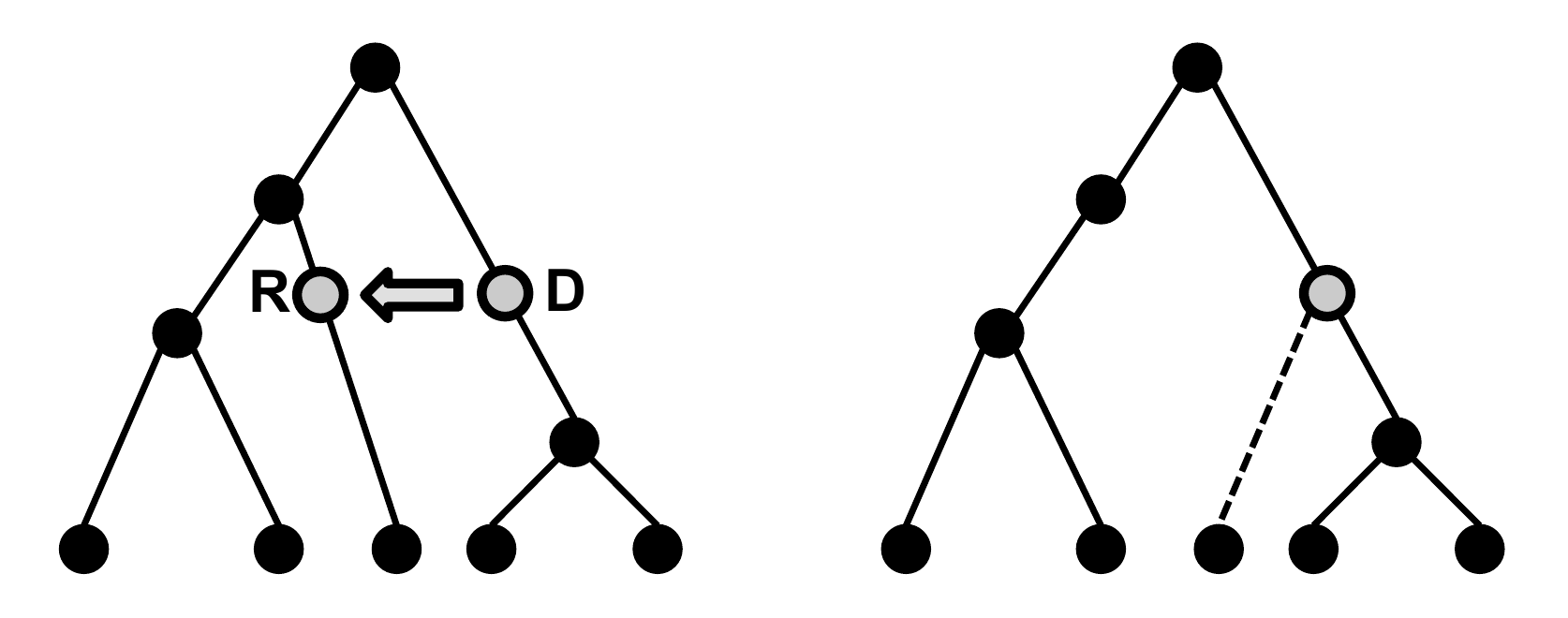}\\
\caption{
An LGT event. On the left, the species
phylogeny is shown with the donor (D) and
recipient (R) locations.
On the right, the resulting (unweighted) gene tree
is shown after the LGT transfer.
}\label{fig:spr}
\end{figure}

Our model of LGT is the following.
Note first that, from a topological point of view, 
an LGT transfer is
equivalent to a subtree-prune-and-regraft (SPR)
operation~\cite{SempleSteel:03}.
The recipient location, that is, the location receiving the
genetic transfer, is the point of pruning.
Similarly, the donor location is the point of
regrafting. 
In other words, 
on the gene tree, a new internal node is created 
at the donor location with two children nodes, 
one being the 
original endpoint of the corresponding edge 
and the other being the node immediately 
under the recipient location in the species phylogeny. 
The original edge going to the latter node is removed.
See Figure~\ref{fig:spr}.
\begin{definition}[Random LGT]
\label{def:randomlgt}
Let $0 < R \leq +\infty$ possibly
depending on $n$ (i.e. not necessarily a constant) and note that we explicitly allow
$R = +\infty$.
Let $T_s = (V_s,E_s,L_s;\root,\time)$
be a fixed species phylogeny.
Let $0 < p \leq 1$ be a sampling effort probability.
A gene tree topology $\tcal_g$ is generated according
to the following continuous-time stochastic process
which gradually modifies the species phylogeny 
starting at the root. There
are two components to the process:
\begin{enumerate}
\item {\bf LGT locations.} 
The recipient and donor locations
of LGT events are selected as follows:
\begin{itemize}
\item {\em Recipient locations.}
Starting from the root,
along each branch $e$ of $T_s$, 
locations are selected as recipient of a genetic tranfer
according to a continuous-time Poisson
process with rate $\lgt(e)$. 
Equivalently, the total number of 
LGT events is Poisson with mean $\lgttotalextinct$
and each such event is located independently
according to the following density.
For a location $x$ on branch $e$, 
the density at $x$ is
$\lgtweight(e)/\lgttotalextinct$.

\item {\em Donor locations.} If $x$ is selected
as a recipient location, the corresponding
donor location $y$ is chosen uniformly at random
in $\ball{x}{R}$. The LGT transfer is then
obtained by performing an SPR move from $x$ to $y$, 
that is, the subtree below $x$ in $T_s$ is moved
to $y$ in $T_g$. 
Note that we perform genetic transfers  
chronologically from the root. 

\end{itemize}

\item {\bf Taxon sampling.} Each extant leaf
is kept independently with probability $p$.
(One could also consider a different
probability for each leaf. We use a fixed
sampling effort $p$ for simplicity.)
The set of leaves selected is denoted by $L_g$.
The final gene tree $T_g$ is then obtained
by keeping the subtree restricted to $L_g$. 
\end{enumerate}
The resulting (random) gene tree topology is denoted by
$\tcal_g$.
\end{definition}
When $R < +\infty$
a transfer can only occur between sufficiently
closely related species. One
could also consider more general donor
location distributions. See e.g.~\cite{PuWoKo:10}.
In Section~\ref{sec:highways}, we consider
a different form of preferential exchange,
highways of gene sharing.

\subsection{Recovering the tree-like trend: Main results}

\paragraph{Problem statement} 
Let $T_s = (V_s,E_s,L_s;\root,\time)$
be an unknown species phylogeny. Using homologous gene sequences for
every gene at hand,
we generate $N$ independent gene tree topologies
$\tcal_{g_1}, \ldots, \tcal_{g_N}$ as above. 
Given the gene trees (or their topologies),
we seek to reconstruct the 
topology $\tcal_s^+ = \tcal[T_s^+]$ 
of the extant phylogeny $T_s^+$.
More precisely
we are interested in 
the amount of LGT that can be sustained
without obscuring the phylogenetic signal.
To derive asymptotic results about this question, 
we make
some assumptions on the underlying phylogeny.
We discuss two cases in detail.

In practice, one estimates
gene trees from sequence data.
We come back to gene tree estimation issues
below.

\paragraph{Bounded-rates model} The following
assumption was introduced in~\cite{DaskalakisRoch:10}
and is related to a common assumption in the
mathematical phylogenetics literature.
\begin{definition}[Bounded-rates model]
\label{def:brm}
Let $0 < \ratiolgt < 1$
and $0 < \ratiotime < 1$
be constants.
Let further
$0 <  \maxtime < +\infty$
be a constant
and 
$0 <  \maxlgt < +\infty$
be a value possibly depending on $n^+$.
Under the Bounded-rates model, we consider
the set of phylogenies 
$T_s = (V_s,E_s,L_s;\root,\time)$ 
with $n^+ > 0$ extant leaves and 
$n^- \geq 0$ extinct leaves
and extant phylogeny 
$T_s^+ = (V_s^+,E_s^+,L_s^+;\root,\time^+)$ 
such that the following conditions are satisfied: 
$$
\minlgt \equiv \ratiolgt \maxlgt \leq \lgt(e) \leq \maxlgt, 
\quad \forall e\in E_s,
$$ 
and 
$$
\mintime \equiv \ratiotime \maxtime \leq \time^+(e^+) \leq \maxtime, 
\quad \forall e^+\in E_s^+.
$$
\end{definition}

Our result in this case is the following.
We use $\maxlgt$ to control the amount of LGT
in the model.
\begin{theorem}[Main result: Bounded-rates model, $R = +\infty$]
\label{thm:main1}
Let $R = +\infty$.
Under the Bounded-rates model,
it is possible to reconstruct 
the topology of the extant phylogeny 
with high probability (w.h.p.)
from $N= \Omega(\log n^+)$ gene tree topologies if
$\maxlgt$ is such that
$$
\lgttotal = O\left(\frac{n^+}{\log n^+}\right).
$$ 
\end{theorem}

In words, we can reconstruct the species
phylogeny w.h.p.~as long as
the expected number of LGT events 
$\lgttotal$ (as measured on the extant phylogeny)
per gene
is at most of the order of $\frac{n^+}{\log n^+}$.
This result is based on a polynomial-time algorithm we describe
in Section~\ref{sec:analysis}.
Note that, in typical phylogenomic studies,
the number of genes is much larger than the number
of species. Therefore, our assumption that the
number of genes should be at least of the
order of the logarithm of the number of
extant species is mild.

We also show that the bound on $\lgttotal$
in Theorem~\ref{thm:main1}
is close to optimal, up to logarithmic factors.
\begin{theorem}[Non-recoverability]
\label{thm:counter}
Under the Bounded-rates model as above
with $N = O(\log n^+)$, 
the topology of the extant phylogeny
cannot, in general,
be reconstructed w.h.p.~if 
$\maxlgt$ is such that 
$\lgttotal = \Omega(n^+ \log \log n^+)$.
\end{theorem}
More generally, the species phylogeny
cannot be reconstructed from $N$ genes if
$\lgttotal = \Omega(n^+ \log N)$.
Theorem~\ref{thm:counter} is proved by a coupling
argument~\cite{Lindvall:92}. 
In words we show that, 
with the order of $\Omega(n^+ \log \log n^+)$
expected LGT events, there is insufficient signal 
from the gene trees 
to distinguish between two species phylogenies
with high probability.

\paragraph{Yule process} Branching processes
are commonly used to model species phylogenies~\cite{RannalaYang:96}. In the
continuous-time Yule process (or pure-birth process), 
one starts with two species (representing
the two branches emanating from the root).
At any given time, each species generates 
a new offspring at rate $0 < \nu < +\infty$. We stop the 
process when the number of species is exactly $n + 1$
(and ignore the $n+1$st species).
This process
generates a species phylogeny with 
$n = n^+$ extant species with branch lengths given
by the inter-speciation times in the above process.
Note that $n^- = 0$ by construction.
Let $0 < \ratiolgt < 1$ be a constant.
We also assume that
$$
\minlgt \equiv \ratiolgt \maxlgt \leq \lgt(e) \leq \maxlgt, 
\quad \forall e\in E_s,
$$ 
for some $0 < \maxlgt < +\infty$ possibly depending
on $n$. As above,
we use $\maxlgt$ to control the amount of LGT
in the model.

An advantage of the Yule model is that,
unlike the Bounded-rates model,
it does not place arbitrary constraints
on the inter-speciation times. In particular,
the following analog of Theorem~\ref{thm:main1} 
suggests that our analysis does not
rely on such constraints.
\begin{theorem}[Main result: Yule process, $R = +\infty$]
\label{thm:main2}
Let $R = +\infty$.
Under the Yule model,
the following holds with probability 
arbitrarily close to $1$.
It is possible to reconstruct 
the topology of the extant phylogeny w.h.p.~from $N = \Omega (\log n)$ gene tree topologies
if 
$\maxlgt$ is such that
$$
\lgttotal =O \left(\frac{n}{\log n}\right).
$$
\end{theorem}

\paragraph{Preferential LGT} 
When $R < +\infty$, that is, when transfers occur
only between sufficiently related species,
we obtain the following generalization which 
implies that preferential LGT makes the
tree-building problem easier.
\begin{theorem}[Preferential LGT]
\label{thm:main3}
Let $0 < R < \log n^+$ possibly depending on $n^+$.
Under the Bounded-rates model, 
it is possible to reconstruct 
the topology of the extant phylogeny w.h.p.~from $N = \Omega(\log n^+)$ gene tree topologies
if 
$\maxlgt$ is such that
$$
\lgttotal = O\left(\frac{n^+}{R}\right).
$$ 
A similar result holds under the Yule model.
\end{theorem}

\paragraph{Further results}
We also obtain results on highways of LGT as well as
sequence-length requirements. 
These results require additional background. See 
Sections~\ref{sec:highways} and~\ref{sec:sequence}
respectively.

\section{Probabilistic Analysis}
\label{sec:analysis}

We assume that we are given $N$ independent gene tree topologies
$\tcal_{g_1}, \ldots, \tcal_{g_N}$ as above. 
Our goal is to reconstruct the extant phylogeny.

Different algorithms are possible.
A simple approach is to take
a majority vote over all gene tree topologies.
But this approach is problematic under taxon
sampling and cannot sustain the high levels
of LGT we consider below. 

Instead we consider
approaches that aggregate partial information
over all gene trees.  
We focus on subtrees over
four taxa whose topologies are called quartets~\cite{SempleSteel:03}.
We show that 
computationally efficient 
quartet-based approaches can sustain
high levels of LGT. Although we prove
our results for the specific method described below,
our analysis is likely to apply to related methods.
In Section~\ref{sec:distance}, we also
give a similar analysis for a distance-based
method of Kim and Salisbury~\cite{KimSalisbury:01}.

\subsection{Algorithm}
\label{sec:algorithm}

We consider the following
approach related to an algorithm of
Zhaxybayeva et al.~\cite{ZhGoCh+:06}. 
Let $X = \{a,b,c,d\}$ be a four-tuple of extant species
The topology $\tcal|X$
of a tree $\tcal$ restricted to $X$
can be summarized with
a {\em quartet split}, or {\em quartet} for short. 
There are three possible (resolved) quartets which
we denote $q_1 = ab|cd$, $q_2 = ac|bd$, and
$q_3 = ad|bc$. We first compute the frequency of
each quartet  over all gene trees displaying $X$, that is,
over all gene trees $g$ such that $X \subseteq L_{g}$,
$$
f_X(q_1) = 
\frac{|\{g_i : X \subseteq L_{g_i},\ \tcal_{g_i}|X = q_1\}|}{|\{g_i : X \subseteq L_{g_i}\}|},
$$ 
and similarly for $q_2, q_3$. 
(We set the frequency to $0$ if the denominator is $0$.)
For each $X$, we choose the quartet  with
highest frequency (breaking ties arbitrarily).
\begin{definition}
A set of quartets $Q=\{q_i\}$,
with $L_{q_i}$ the leaf set of $q_i$, 
is {\em compatible} 
if there is a tree ${\cal T}$ 
with leaf set $L_Q \equiv 
\cup_{q_i \in Q}
L_{q_i}$ such that ${\cal T}$ agrees
with every $q_i$.
\end{definition}
Quartet compatibility is, in general, NP-hard~\cite{NPsteel}. However,
when the set $Q$  covers
all possible four-tuple of taxa 
(that is, exactly $n \choose 4$ quartets
with no repeated four-tuple of taxa), there is
a polynomial-time algorithm 
for compatibility~\cite{BanDre,Buneman1971,Berr2001}.
In our procedure, 
for every four-tuple of taxa, there is a single quartet
chosen, so we can
check compatibility easily and output the corresponding
tree. 
In practice, if $Q$ is not compatible, one can use instead
a heuristic supertree method such as 
MRP~\cite{Ragan1992,Baum1992}
or Quartet
MaxCut~\cite{Snir-TCBB-2010,Snir-MPE-2011}. 

The algorithm, which we call QuartetPlurality (QP), is detailed in Figure~\ref{fig:plurality}.
\begin{figure*}[!ht]
\framebox{
\begin{minipage}{13cm}
{\small \textbf{Algorithm} QuartetPlurality\\
\textit{Input:} Gene trees $g_1, \ldots, g_N$;\\
\textit{Output:} Estimated species phylogeny $\hat{T}$;

\begin{itemize}
\item Set $Q = \emptyset$
\item For all four-tuple of taxa $X = \{a,b,c,d\}$, letting
$q_1 = ab|cd$, compute
$$
f_X(q_1) = 
\frac{|\{g_i : X \subseteq L_{g_i},\ \tcal_{g_i}|X = q_1\}|}{|\{g_i : X \subseteq L_{g_i}\}|},
$$ 
and similarly for  $q_2 = ac|bd$ and
$q_3 = ad|bc$. Add the quartet  with
highest frequency (breaking ties arbitrarily)
to $Q$.

\item Using Buneman's algorithm~\cite{Buneman1971}
compute the tree $\hat{T}$ compatible
with $Q$ (or abort if no such tree is found).

\item Output $\hat{T}$.

\end{itemize}

}
\end{minipage}
} \caption{
Algorithm QuartetPlurality.
} \label{fig:plurality}
\end{figure*}

\subsection{A general formula}

Our asymptotic analysis is based on the following
claim.  Recall that, for a subset of extant species $X$, 
we let
$\tcal_s|X$ be the extant phylogeny topology restricted 
to $X$ with corresponding edge set $\ecal(\tcal_s|X)$. Also recall
that  $\lgtweight(e) =  \lgt(e) \time(e)$ is the expected number of
LGT events on \Red{edge} $e$ which we refer to as the 
{\em LGT weight}, or {\em weight} for short, of $e$.
Let
$$
\lgttotal_X = \sum_{e \in  \ecal(\tcal_s|X)} \lgtweight(e),
$$ 
be the total weight of the subtree $\tcal_s|X$
under the weights $\lgtweight(e)$, $e\in E_s$.
Define the {\em maximum quartet weight (MQW)} as
$$
\mrqw 
= \max\{\lgttotal_X : X \subseteq (L_s^+)^4\}.
$$
\begin{lemma}[Probability of a miss]
\label{lem-missq}
Let $\tcal_g$ be a gene tree topology 
distributed according 
to the random LGT model such that
$X = \{a,b,c,d\} \subseteq L_g$. 
Let $q^X_s$ (respectively $q^X_g$) be the 
quartet corresponding to $\tcal_g|X$ 
(respectively $\tcal_s|X$). Then
$$
\P[q^X_g = q^X_s | X \subseteq L_g]
\geq \exp\left(-\mrqw\right).
$$
\end{lemma}
Recall that $\lgttotal$ is the expected number of
LGT events (as measured on the extant phylogeny) 
per gene.
As a comparison, 
note that the probability that a gene tree
is LGT-free is $e^{-\lgttotal}$,
which can be much smaller. 

\begin{prevproof}{Lemma}{lem-missq}
We first note that,
by our assumption
that the species phylogeny is bifurcating,
$q^X_s$ is resolved.
Similarly $q^X_g$ is resolved because
under a Poisson process for the
recipient location the probability that
a vertex has degree higher than $2$
(that is, that a pruning and re-grafting occurs exactly
at the location of an existing vertex)
is $0$. 

Now we observe that if none of the recipient
locations lands on $\tcal_s|X$ then the corresponding
quartet remains intact. Indeed an SPR
move can only (potentially) affect those quartets with at least 
one leaf in the pruned subtree\Red{, and this happens with probability
 $ \frac{\lgttotal_X}{\lgttotal}$}. 
The claim
then follows by induction on the number
of LGT events. 

Hence the probability that $q^X_g = q^X_s$
is at least the probability that all LGT events
(on the extant phylogeny)
miss $\tcal_s|X$, which is at least
\begin{eqnarray*}
\P[q^X_g = q^X_s | X \subseteq L_g]
&\geq& \sum_{i=0}^{+\infty} \frac{e^{-\lgttotal} \lgttotal^i}{i!}
\left(1 - \frac{\lgttotal_X}{\lgttotal}\right)^i\\
&=& e^{-\lgttotal} 
\exp\left(\lgttotal\left(1 - \frac{\lgttotal_X}{\lgttotal}\right)\right)\\
&\geq& \exp\left(-\mrqw\right).
\end{eqnarray*}
\end{prevproof}

\subsection{Bounded-rates and Yule models}

Next we argue that, 
under appropriate assumptions on the species phylogeny, 
the
maximum quartet weight is bounded in such a way 
that the plurality
quartet topology for every 
four-tuple of taxa $X = \{a,b,c,d\}$, 
which we denote by $q^X_*$,
satisfies $q^X_* = q^X_s$. 
As a result, our quartet set is
compatible and $\tcal_s^+$ can be reconstructed 
efficiently.

\subsubsection{Bounded-rates model}
We bound the maximum quartet weight $\mrqw$ in the 
Bounded-rates model.
\begin{lemma}[Bound on quartet weight: Bounded-rates case]
\label{lem-relbrm}
Under the Bounded-rates model 
it holds that
$$
\mrqw = O\left(\maxlgt \log n^+\right),
\qquad \lgttotal = \Theta(\maxlgt n^+).
$$
\end{lemma}
\begin{prevproof}{Lemma}{lem-relbrm}
The first part of the proof is taken
from~\cite{DaskalakisRoch:10}.
Let $h$ (respectively $H$) be the smallest 
(respectively largest)
number of edges on a
path between the root and an extant leaf.
Because the number of extant leaves is $n^+$ 
and the extant phylogeny is bifurcating
(recall that we suppressed vertices of degree $2$
after taking a restriction to the extant species),
we must have
$2^h \leq n^+$ and $2^H \geq n^+$. 
Since all extant leaves are contemporaneous 
it must be that 
$H \mintime \leq h \maxtime$. Combining these constraints gives
\begin{equation*}
\frac{\mintime}{\maxtime} \log_2 n^+ 
\leq h 
\leq H 
\leq \frac{\maxtime}{\mintime} \log_2 n^+.
\end{equation*}
Hence 
$$
\max\{\lgttotal_X : X \subseteq (L_s^+)^4\}
\leq 4 \maxlgt \maxtime \frac{\maxtime}{\mintime} \log_2 n^+.
$$

The total number of edges 
in the extant phylogeny
is $2n^+ - 3$ so that
$$
\lgttotal = \Theta(\maxlgt n^+).
$$
\end{prevproof}

Using Lemma~\ref{lem-relbrm}, we prove
Theorem~\ref{thm:main1}. 
First recall the following standard concentration inequality (see e.g.~\cite{MotwaniRaghavan:95}):
\begin{lemma}[Azuma-Hoeffding Inequality]\label{lemma:azuma}
Suppose ${\bf Z}=(Z_1,\ldots,Z_m)$ 
are independent random variables taking values in a set
$S$, and $h:S^m \to \real$ is any $t$-Lipschitz function:
$|h(\mathbf{z}) - h(\mathbf{z'})|\leq t$ 
whenever $\mathbf{z}, \mathbf{z'} \in S^m$ 
differ at just one coordinate. Then,
$\forall \zeta > 0$,
\begin{equation*}
\prob\left[|h({\bf Z}) - \expec[h({\bf Z})]| \geq \zeta \right]
\leq 2\exp\left(-\frac{\zeta^2 }{2 t^2 m}\right).
\end{equation*}
\end{lemma}

\begin{prevproof}{Theorem}{thm:main1}
Consider the quartet-based approach described in
Section~\ref{sec:algorithm}. 
Take $\maxlgt = C_1/\log n^+$ 
with $C_1 > 0$ small enough so that
$$
\lgttotal =O\left(\frac{n^+}{\log n^+}\right),
$$ 
and using Lemmas~\ref{lem-missq} and~\ref{lem-relbrm},
we have for any four-tuple $X$ of extant species
$$
\P[X \subseteq L_g]
= p^4,
$$
and
$$
\P[q^X_g = q^X_s\,|\,X \subseteq L_g]
\geq \exp\left(-\mrqw\right)
\geq \exp\left(-O(C_1)\right)
\geq \frac{2}{3},
$$
for $C_1$ small enough.
We choose $C_2 > 0$ large enough with
$$
N \geq C_2 \log n^+,
$$
and $\eps < p^4$
so that, using Lemma~\ref{lemma:azuma},
the following inequalities hold.
Consider the following events
\begin{eqnarray*}
\ecal_0 = \{
||\{g_i : X \subseteq L_{g_i}\}| - N p^4|\leq N \eps\}
\end{eqnarray*} 
and
\begin{eqnarray*}
\ecal_1 = \left\{
|\{g_i : X \subseteq L_{g_i},\ \tcal_{g_i}|X = q_1\}| > 
\frac{1}{2} |\{g_i : X \subseteq L_{g_i}\}|\right \}.
\end{eqnarray*} 
By Lemma~\ref{lemma:azuma},
$$
\P[\ecal_0^c]
\leq  \exp\left(- O(\eps^2 N) \right),
$$
and
$$
\P[\ecal_1^c\,|\,\ecal_0]
\leq \exp\left(-O(N(p^4 - \eps))\right).
$$
Hence, for a constant $C_2$ large enough,
\begin{eqnarray*}
\P[f_X(q^X_s) < 1/2] 
&\leq& \P[\ecal_0^c] + \P[\ecal_1^c\,|\,\ecal_0]\\
&\leq& O\left(\frac{1}{(n^+)^4}\right).
\end{eqnarray*}
Then the plurality vote is correct for every four-tuple
of taxa
and the extant phylogeny is correctly reconstructed.
\end{prevproof}

\subsubsection{Yule process}
We now consider the 
Yule model.
\begin{lemma}[Bound on quartet weight: Yule case]
\label{lem-relyule}	
Under the Yule model, it holds that
$$
\mrqw = \Theta\left(\maxlgt \log n\right),
\qquad \lgttotal = \Theta\left(\maxlgt n\right)
$$
with probability approaching $1$ as $n \to +\infty$.
\end{lemma}
\begin{prevproof}{Lemma}{lem-relyule}
We consider a pure-birth process with 
birth rate $\nu$ starting from $2$ species.
For background on branching processes see~\cite{AtheryaNey:72}.

Let $Z_i$ be the $(i-1)$-th inter-speciation time.
As a minimum of 
$i$ independent exponential distributions
with mean $1/\nu$,
$Z_i$ is an exponential with
mean $1/(i\nu)$.
Moreover the $Z_i$s are independent.
Hence the height of the phylogeny in time units, that is, 
the total time until $n+1$ species
are present (recall that we ignore the $(n+1)$-st
species) is
$$
{\bf Z}
= \sum_{i=2}^{n+1} Z_i,
$$
and we have
$$
\E[{\bf Z}]
= \sum_{i=2}^{n+1} \E[Z_i]
= \sum_{i=2}^{n+1} \frac{1}{i \nu}
= \Theta(\nu^{-1} \log n),
$$
and
$$
\var[{\bf Z}]
= \sum_{i=2}^{n+1} \var[Z_i]
= \sum_{i=2}^{n+1} \frac{1}{i^2 \nu^2}
= \Theta(\nu^{-2}).
$$
The total weight of the phylogeny
in time units
$$
{\bf Y}
= \sum_{i=2}^{n+1} i Z_i,
$$
is a sum of $n$ independent
exponential random variables with
parameter $\nu$,
and we have
$$
\E[{\bf Y}]
= \sum_{i=2}^{n+1} i \E[Z_i]
= \sum_{i=2}^{n+1} i \frac{1}{i \nu}
= \nu^{-1} n,
$$
and
$$
\var[{\bf Y}]
= \sum_{i=2}^{n+1} i^2 \var[Z_i]
= \sum_{i=2}^{n+1} i^2 \frac{1}{i^2 \nu^2}
= \nu^{-2} n.
$$

By Chebyshev's inequality, 
$$
\P[{\bf Z} \geq C_1 \log n]
\leq \frac{C_2}{C_3 \log^2 n} \to 0,
$$
and
$$
\P[{\bf Y} \leq C_4 n]
\leq \frac{C_5 n}{C_6 n^2} \to 0,
$$
for appropriately chosen $C$s not depending
on $n$.
The same holds in the other direction
so that $\mrqw = \Theta(\maxlgt \log n)$
and $\lgttotal = \Theta(\maxlgt n)$ with
probability approaching $1$.
\end{prevproof}


\begin{prevproof}{Theorem}{thm:main2}
Using Lemma~\ref{lem-relyule},
the proof of Theorem~\ref{thm:main2} follows
form the same lines as that of Theorem~\ref{thm:main1}.
\end{prevproof}

\subsection{Preferential LGT}

We now prove Theorem~\ref{thm:main3}.

\begin{prevproof}{Theorem}{thm:main3}
The proof is similar to that of Theorems~\ref{thm:main1}
and~\ref{thm:main2}. The main difference
is in the proof of Lemma~\ref{lem-missq}. 
In that proof, note that if $R < +\infty$
then for an LGT to affect the quartet on $X$, 
it must
be that not only 1) the recipient location lands on
$\tcal_s|X$, but also 2) that it lands on a location
below either branchings of the 
corresponding quartet tree 
within time $R$ of the branching point. Indeed 
these are the only locations where the 
corresponding leg of the quartet tree
can potentially jump to a subtree corresponding to a
different leg. 
(In fact, it must be that a leg {\em on the other side
of the internal branch of the quartet tree} is within
time $2R$.)
The length of this region 
is at most $4R$ in $\time$-distance.
Hence in the bound on the probability of a miss we
get
\begin{eqnarray*}
\P[q^X_g = q^X_s  | X \subseteq L_g]
&\geq& \exp\left(-\min\{\mrqw,4 R \maxlgt \}\right).
\end{eqnarray*}
The result then follows.
\end{prevproof}

\subsection{Non-recoverability}

We now prove Theorem~\ref{thm:counter}.

\begin{prevproof}{Theorem}{thm:counter}
We use a coupling argument~\cite{Lindvall:92}.
Fix $\delta > 0$ small.
We construct two species phylogenies with
different topologies which cannot be
distinguished with probability $1-\delta$ 
from $N$ gene tree topologies
when the total expected
amount of LGT $\lgttotal$ 
is of the order of $n^+ \log\log n^+$
per gene. In particular the reconstruction
problem cannot be solved in that case.
The idea of a coupling is to
run the stochastic processes of LGT on 
both phylogenies simultaneously so as to output
the same gene trees with high probability
without changing the marginal distributions (that is,
the probability distributions of gene tree topologies
on each phylogeny separately).

We proceed as follows. Consider a
complete binary tree $T_s'$ 
on a set of $n$ leaves (all extant) and denote the four children at
height 2 from the root as $a, b, c,d$, where $a$ and $b$ are sisters and so are $c$
and $d$. Let $T_z$ be the subtree with $n/4$ leaves rooted at $z\in \{a,b,c,d\}$. Moreover, for simplicity, assume all edges of
$T_s'$ have the same LGT weight. From  $T_s'$ we construct $T_s''$
by rewiring the four nodes $\{a,b,c,d\}$ such that $a$ is now
sister with $c$ and $b$ with $d$.

We generate 
$N = \Theta(\log n)$ genes \Red{trees} on \Red{each of $T_s'$ and $T_s''$}
as follows. We run the stochastic 
process of LGT on $T_s'$ as described 
in Definition~\ref{def:randomlgt}. 
Let
$\tcal_{g_1}',\ldots,\tcal_{g_N}'$
be the gene tree topologies so obtained.
For $T_s''$ and every gene,
we use {\em exactly the same LGT events} as the ones
generated on $T_s'$ where we identify the
two edges adjacent to the roots in $T_s'$ and
$T_s''$ arbitrarily. 
Let
$\tcal_{g_1}'',\ldots,\tcal_{g_N}''$
be the gene tree topologies so obtained.

Since $T_s'$ and $T_s''$ are identical
below every 
$z\in \{a,b,c,d\}$ and LGT events occur only between
contemporaneous points, the subtrees under $\{a,b,c,d\}$ 
in $\tcal_{g_i}'$ and
$\tcal_{g_i}''$ are identical for every gene $i$.

For $z \in Z$, let $e_z$ be the edge adjacent to $z$
and above it in $T_s'$ (and in $T_s''$). It remains to show that, for $\tcal_{g_i}'$ and $\tcal_{g_i}''$
to be identical under the joint construction
above, it suffices that the following
{\em good event} occurs: three consecutive 
LGT moves start on the {\em same} edge 
in $e_a, \ldots, e_d$ (donor location)
and land on the other three edges 
in $e_a, \ldots, e_d$ (recipient location), for example,
$a\rightarrow d$, $a\rightarrow c$, $a\rightarrow b$.
See Figure~\ref{fig:nonrecover}. 
\begin{figure}[t]
\centering
\includegraphics[width=4in]{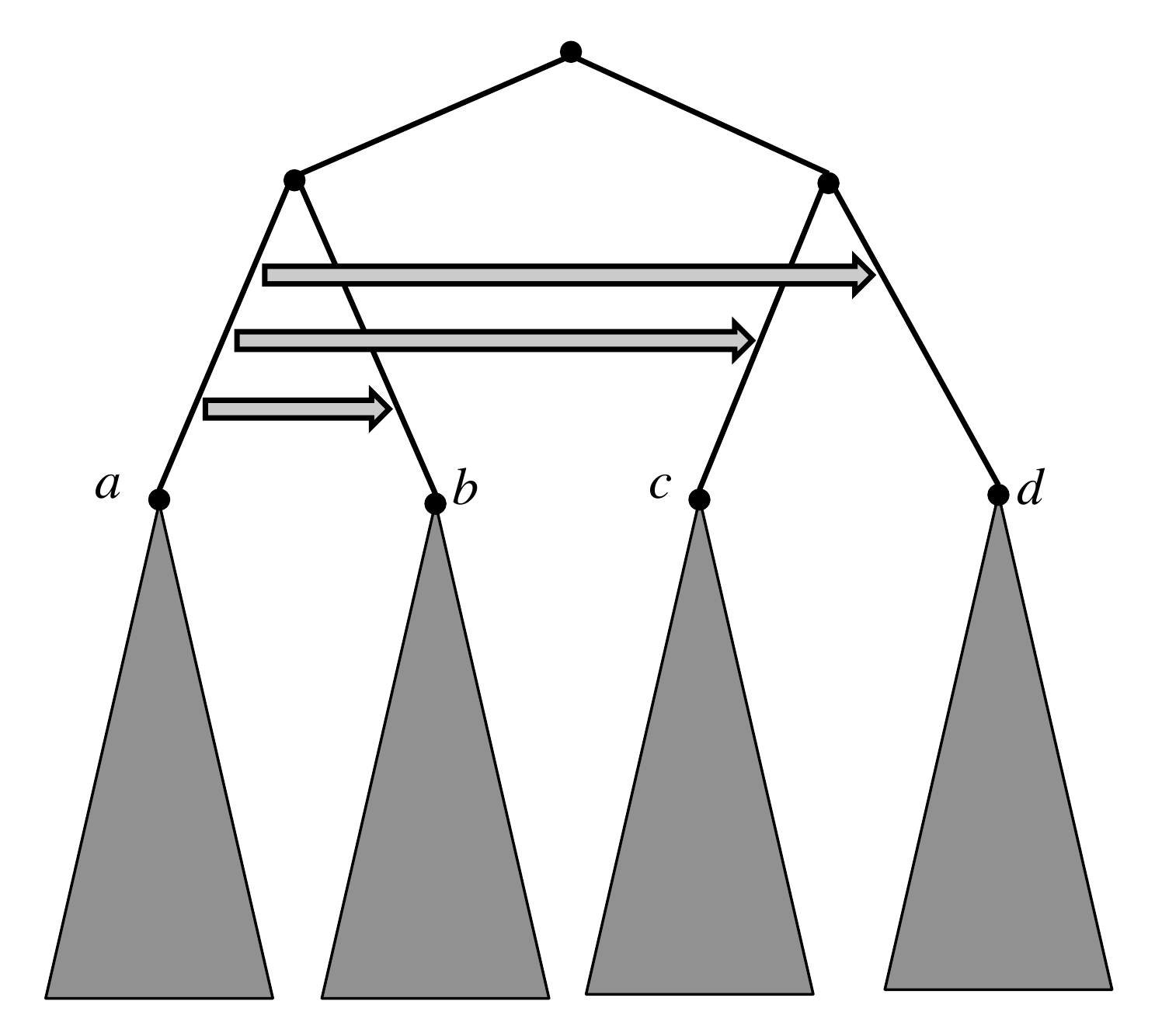}\\
\caption{
Good event.
}\label{fig:nonrecover}
\end{figure}
Indeed, in that case, 
the first donor location above becomes the common 
ancestor to all nodes in the gene trees. 
From that point on, we obtain
the same gene tree for both phylogenies.

We claim that the probability that the good event
does not occur is $O(1/\log n)$.
Under the assumption that $\lgttotal = \Omega(n \log \log n)$
and that the LGT weights are equal,
the number of LGT events on any edge
is Poisson with mean $\Omega(\log \log n)$.
Consider the time interval between
the nodes at height $1$ from the root 
and the nodes at height $2$.
Divide this interval into $\nu = O(\log \log n)$
equal subintervals $I_1, \ldots, I_\nu$ 
such that the number of
LGT events on edge $e_z$ in $I_i$ is Poisson
with mean $C_0$ for some constant $C_0 > 0$. 
In $I_i$ the probability
that there is no LGT event originating from $e_b,
\ldots, e_d$ and that there is exactly three
LGT events originating from $e_a$ 
and landing on $e_b, e_c, e_d$ in that order is
$$
\tilde{p}
= \left(e^{-C_0}\right)^3
\left(e^{-C_0} \frac{C_0^3}{ 3!} \left(\frac{1}{3}\right)^3\right)
\equiv C_1.
$$
The subintervals are independent.
The probability
that the event above does not happen
in any of $I_1,\ldots,I_\nu$, is 
$$
\tilde{p}^\nu = (1-C_1)^\nu = O\left(\frac{1}{\log n}\right).
$$
This gives an upper bound of $O(1/\log n)$
on the probability that the good event does not happen.

Therefore, by a union bound over the genes, the probability that the good event
does not occur on at least one gene tree
is $\Theta(\log n)\cdot O(1/\log n)=O(1)$, which
is at most
$\delta$ if the constant in
$\lgttotal$
is large enough. 
If the good event occurs on every gene tree, 
then both phylogenies
output the exact same set of gene tree topologies.
That concludes the proof.
\end{prevproof}

\section{Highways of LGT}
\label{sec:highways}

In this section, we add highways of gene sharing
to the model. Highways are, in essence,
non-random patterns of LGT~\cite{BeHaRa:05}. 
These
can potentially take different shapes. 
Following Bansal et al.~\cite{BaBaGo+:11},
we focus on pairs of edges in the 
phylogeny that undergo an unusually
large number of LGT
events between them. 

We give two results. As long as
the frequency of genes affected by highways is low 
enough, the species phylogeny can be reconstructed
using the same approach as in Section~\ref{sec:analysis}.
Moreover, with extra assumptions on the positions
of the highways with respect to each other, 
the highways themselves can be inferred.

In this section, we assume $n^- = 0$.

\subsection{Model}

We generalize our model of LGT as follows.
\begin{definition}[Highways of LGT]
Let $T_s = (V_s,E_s,L_s;\root,\time)$ be a species
phylogeny with LGT rates $0 < \lgt(e) < +\infty$, $e\in E_s$
and let $0 < p \leq 1$ be a taxon sampling probability.
Assume $n^- = 0$.
For $\beta = 1,\ldots, B$, let 
$\hw_\beta = (e^H_{\beta,0}, e^H_{\beta,1})$ 
be a pair of
edges in $T_s$ which share contemporaneous
locations.
We call $\hw_\beta$ a highway. 
Let $g_1, \ldots, g_N$
be $N$ genes. 
Each highway $\hw_\beta$ involves 
a subset $\hwgenes_\beta$ of
the genes. 
If gene $g_i \in \hwgenes_\beta$, then it 
undergoes an LGT event between 
a pair of contemporaneous 
locations $x^H_{\beta,i} \in e^H_{\beta,0}$
and $y^H_{\beta,i} \in e^H_{\beta,1}$.
We let $\hwfrac_\beta$ be the fraction
of genes such that $g_i \in \hwgenes_\beta$
and we assume that $\hwfrac_\beta > \minhwfrac$
for some $\minhwfrac$ (chosen below).
In addition, independently
from the above, we assume that each gene
undergoes LGT events at random locations
as described in Definition~\ref{def:randomlgt}.
We denote by $\tcal_{g_1},\ldots,\tcal_{g_N}$
the gene tree topologies so obtained.
\end{definition}
\begin{remark}[Deterministic setting]
Note that the highways and which genes
are involved in them are deterministic in this setting.
Only the random LGT events are governed by
a stochastic process.
Note moreover that we allow highway events
to go in either direction, that is, from $e^H_{\beta,0}$
to $e^H_{\beta,1}$ or vice versa.
\end{remark}

\subsection{Building the species tree in the
presence of highways}

We first prove that the species phylogeny
can still be reconstructed in the presence
of highways as long as the fraction of genes
involved in highways is low enough. 
We only discuss the Bounded-rates
model with $R = +\infty$. 
\begin{theorem}[Highways of LGT]
\label{thm:main4}
Consider the Bounded-rates model
with $R = +\infty$ and assume that $B < +\infty$
is constant. Assume further that there 
is a constant $0 < \maxhwfrac < 1$
such that
$$
\hwfrac_\beta < \maxhwfrac, 
\quad 
\beta=1,\ldots,B.
$$
If 
$$
\maxhwfrac < \frac{1}{2B},
$$ 
then it is possible to reconstruct 
the topology of the extant phylogeny w.h.p.~from $N  =\Omega (\log n^+)$ gene tree topologies
if 
$\maxlgt$ is such that
$$
\lgttotal = O\left(\frac{n^+}{\log n^+}\right).
$$ 
\end{theorem}
\begin{prevproof}{Theorem}{thm:main4}
The proof is similar to that of Theorem~\ref{thm:main1}.
Note that a quartet tree in the
species phylogeny can be affected by a highway
in at most a fraction $< B \frac{1}{2B} = \frac{1}{2}$
of the genes. Moreover by the
proof of Lemma~\ref{lem-missq}, 
choosing $C_1$ small
enough, a quartet tree is affected by a random
LGT event in an arbitrarily small fraction of genes.
Therefore the plurality vote will reconstruct
the correct split with high probability. The result
follows. 
\end{prevproof}

\subsection{Inferring highways}

The problem of inferring the highway 
locations is essentially a network reconstruction
problem. Such problems are often
computationally intractable. See e.g.~\cite{HuRuSc:10}.
Therefore, we require some extra assumptions.
Our goal here is not to provide the most general
result but rather to illustrate that our analysis
extends naturally to certain network settings.
The following assumption is related to so-called
galled trees.
\begin{assumption}
\label{assump:galled}
We assume that no highway connects
two edges in $T_s$ separated by less than two
edges
or edges adjacent to root edges.
(Such cases cannot be reconstructed.)
Seen as an edge superimposed on 
$T_s$,
a highway event $(x^H_{\beta,i}, y^H_{\beta,i})$
forms a cycle. We assume that all such
cycles are disjoint, that is, they do not
share common locations.
\end{assumption}
We then prove the following.
We use a computationally
efficient algorithm, which we call
RoadRoller, described
in Figure~\ref{fig:roller} and explained in the proof.
\begin{figure*}[!ht]
\framebox{
\begin{minipage}{13cm}
{\small \textbf{Algorithm} RoadRoller\\
\textit{Input:} Gene trees $g_1, \ldots, g_N$;\\
\textit{Output:} Estimated species phylogeny $\hat{T}$
and highway locations;

\begin{itemize}
\item Use QuartetPlurality to reconstruct the
species phylogeny $\hat{T}$. Let $\qcal$ be the set of all quartets 
whose estimated frequency is less than $1/2$
but more than $\minhwfrac/2$. 

\item For all pairs of four-tuples $X \neq X'$ (possibly
sharing taxa) with a corresponding
quartet
in $\qcal$, 
\begin{itemize}
\item Find the shared edges $e(X,X')$ 
along the internal branches 
of $\tcal_s|X$ and $\tcal_s|X'$. 
\item Let $X\sim X'$ if
$e(X,X') \neq \emptyset$. 
\end{itemize}
 \item Build the  
graph $\gcal$ 
corresponding to $\sim$
with vertex set being all $X$s with
a corresponding quartet in $\qcal$.

\item For each connected component $W$ of
$\gcal$,
\begin{itemize}
\item Compute the union $\pcal$ of all $e(X,X')$ over
pairs $X$ and $X'$ in $W$. Abort if $\pcal$ is not
a path.
\item Let $\tilde{e}_0^{W}$
and $\tilde{e}_1^{W}$ be the start and
end edges on the path $\pcal$.
\item For $i=0,1$, let $e_i^-$ and
$e_i^+$ be the edges adjacent to $\tilde{e}_i^W$.
\item For each pair with one element in 
$\{e_0^-, e_0^+\}$ and one element in
$\{e_1^-, e_1^+\}$, determine whether 
each $\tcal_s|X$ with $X$ in $W$ 
contains at least one element in the pair.
\item If only one pair passed the previous 
test, 
\begin{itemize}
\item Denote the pair by $(e_0^W,e_1^W)$,
\item Else, let $e_0^W$ be the intersection
of the pairs found (abort if the intersection
does not contain a unique element),
choose an $X$ in $W$ such that $\tcal_s|X$
includes all of $\{e_0^-, e_0^+\}$ and 
$\{e_1^-, e_1^+\}$, and use the corresponding
quartet in $\qcal$ to determine the 
sister leaf to the leaf below $e_0^W$.
The latter leaf is below edge $e_1^W$
among $\{e_0^-, e_0^+,e_1^-, e_1^+\}$.
\end{itemize}
\end{itemize}

\item Output $\hat{T}$ and $(e_0^W,e_1^W)$
for all $W$.
\end{itemize}
}
\end{minipage}
} \caption{
Algorithm RoadRoller.
} \label{fig:roller}
\end{figure*}
\begin{theorem}[Inferring highways]
\label{thm:main5}
Consider the Bounded-rates model
with $R = +\infty$ and assume that $B < +\infty$
is constant. Assume further that there 
are constants $0 < \minhwfrac < \maxhwfrac < +\infty$
such that
$$
\minhwfrac < \hwfrac_\beta < \maxhwfrac, 
\quad 
\beta=1,\ldots,B.
$$
If 
$$
\maxhwfrac < \frac{1}{2},
$$ 
and Assumption~\ref{assump:galled} holds
then 
it is possible to reconstruct
the topology of the extant phylogeny
as well as the highway edges w.h.p.~from $N = \Omega(\log n^+)$ gene tree topologies
if 
$\maxlgt$ is such that
$$
\lgttotal = O\left(\frac{n^+}{\log n^+}\right).
$$ 
\end{theorem}
\begin{prevproof}{Theorem}{thm:main5}
Consider a four-tuple $X$ such that
$\tcal_s|X$ contains at least one
highway location and such that the quartet
$q_s^X$ is modified by the corresponding
highway. Because such a highway must
connect a leg of $\tcal_s|X$ to a subtree
on the other side of the internal branch
of $\tcal_s|X$, our galled tree assumption
implies that any given quartet tree can be
affected by at most one highway,
otherwise the corresponding cycles
would intersect along the internal branch. 
Hence, from the proof of 
Theorem~\ref{thm:main4} and the assumption
that $\maxhwfrac < \frac{1}{2}$
(instead of $\maxhwfrac < \frac{1}{2B}$),
we can
reconstruct the extant phylogeny.

Further, it follows by the proof of
Theorem~\ref{thm:main4} that,
if $\minhwfrac > 0$ and $C_1$ is small enough,
the second most frequent quartet over a four-tuple
as above
is the one obtained by going through the highway.
Let $\qcal$ be the set of all quartets 
whose estimated frequency is less than $1/2$
but more than $\minhwfrac/2$. 
By the previous argument and 
Lemma~\ref{lemma:azuma} (see the 
proof of Theorem~\ref{thm:main1}
for a similar computation),  $\qcal$ contains 
w.h.p.~exactly those
quartets affected by a highway.

For $X, X'$ with quartets in $\qcal$, write $X \sim X'$ if 
the quartet trees $\tcal_s|X$ and
$\tcal_s|X'$ share an edge {\em along
their internal branch}. 
Let $e(X,X')$ be the
set of all such shared edges. 
Note that, although we are considering
four-tuples affected by highways, we are
working on the species phylogeny $\tcal_s$
which has been reconstructed.

By the argument
above, quartets sharing an edge along
their internal branch are necessarily affected
by the same highway. 
Take the 
transitive closure $\sim_*$ of $\sim$. 
Let $W$ be
an equivalence class of $\sim_*$. 
We reconstruct the corresponding highway
as follows. The union of all edges in $e(X,X')$
for some pair $X, X'$ in $W$ forms
a path $\pcal$ in $\tcal_s$. Let $\tilde{e}_0^{W}$
and $\tilde{e}_1^{W}$ be the start and
end edges on this path. The highway corresponding
to $W$ connects an edge $e_0^W$ adjacent to 
$\tilde{e}_0^{W}$
with an edge $e_1^W$ adjacent to 
$\tilde{e}_1^{W}$.
See Figure~\ref{fig:highway}. 
(Note that a highway is represented by
exactly one $W$ because w.h.p.~all
quartets affected by this highway are in $\qcal$
and they are all connected under $\sim$.
See Figure~\ref{fig:highway}.)
\begin{figure}[t]
\centering
\includegraphics[width=4in]{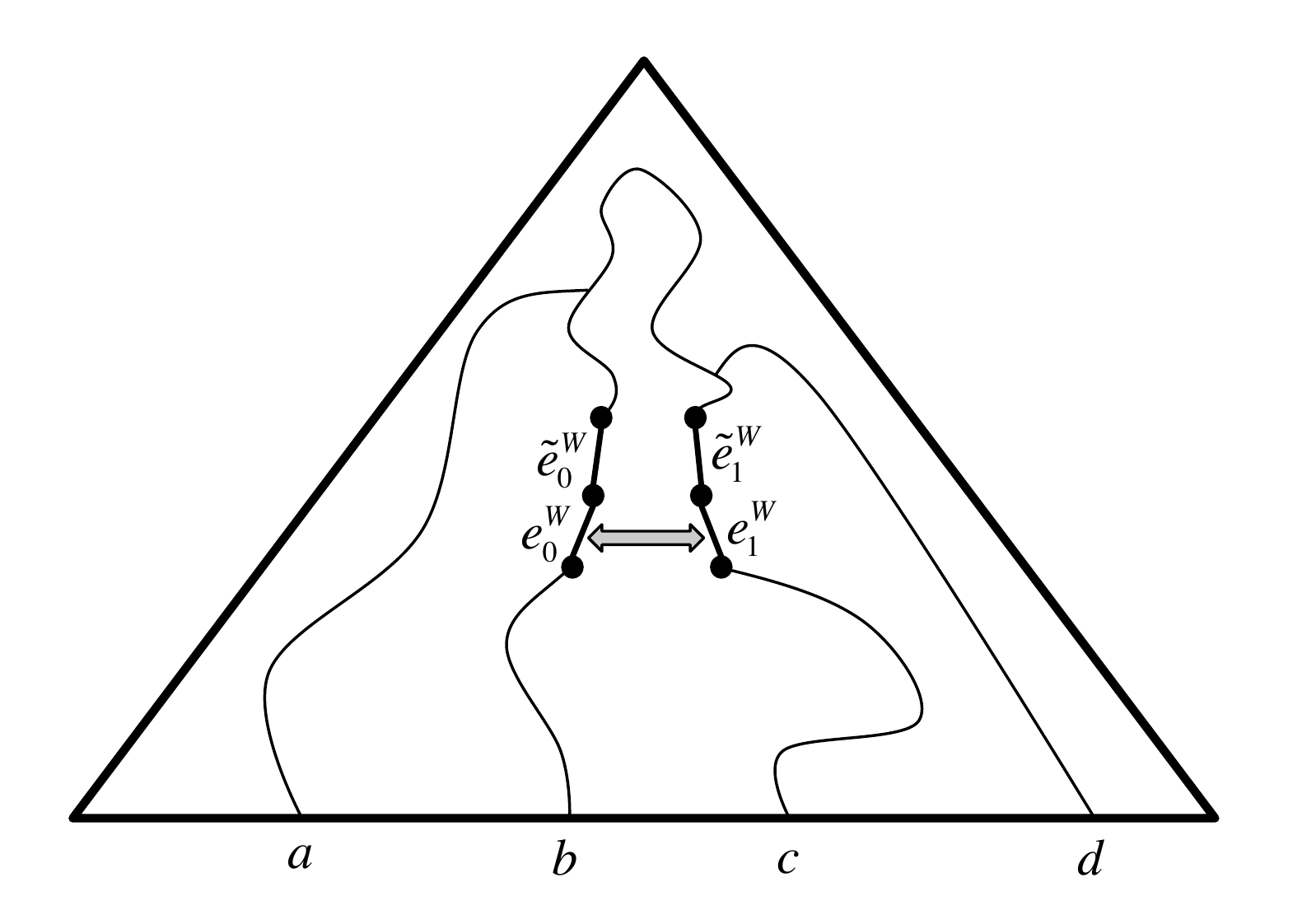}\\
\caption{
Setup in the proof of
Theorem~\ref{thm:main5}. 
The grey arrow
indicates a highway.
Here
$X = \{a,b,c,d\}$, $\tcal_s|X = ab|cd$
and $bc|ad \in \qcal$.
}\label{fig:highway}
\end{figure}

As we argued in the proof of Lemma~\ref{lem-missq},
all quartets affected by the highway corresponding 
to $W$ contain at least one leaf
in a pruned subtree. Because we allow LGT events
in both direction along a highway, there are
two potential pruned subtrees. Moreover, the other
three leaves must be in separate subtrees
hanging from the path $\pcal$. By our assumption,
there are at least three such subtrees (in addition
to the two potentially pruned subtrees). 

Hence, the pruned
subtrees can be identified by checking the
four-tuples in $W$ and finding the pairs of subtrees
with at least one of them present in all of $W$. 
If there is a unique such pair,
this gives the two highway edges and we are done.
Otherwise, the recipient edge is the intersection
of the pairs found.
To identify the donor edge, one simply needs to
use a four-tuple $X$ of leaves in the four
adjacent subtrees to the endpoints of $\pcal$
and check to which branch of $\tcal_s|X$ the subtree corresponding
to the recipient edge is moved in $\qcal$
(that is, in the highway-affected quartet topology).
\end{prevproof}

\section{Distance method and sequence lengths}
\label{sec:sequence}

In this section, in the highway-free case, 
we analyze an alternative, 
distance-based approach that has been considered
in the literature and we provide sequence-length requirements.
Although the quartet-based method analyzed
in Section~\ref{sec:analysis} can in 
principle handle arbitrary branch lengths
(as only the topology of the gene trees is used),
here we need to assume that the gene tree branch lengths
are determined by inter-speciation times and
lineage-specific rates of substitution.
For simplicity, we
assume that there is no gene-specific substitution
rate. 
In practice, one could incorporate such rates
by using a normalization
procedure as detailed 
in~\cite{KimSalisbury:01,GeWaKi:05}.

\subsection{A distance-based approach}
\label{sec:distance}

We analyze a distance-based
approach similar to that introduced 
in~\cite{KimSalisbury:01} and studied
empirically in~\cite{GeWaKi:05}.
Given branch lengths, 
a gene tree is naturally equipped with a
{\em tree metric} on the leaves
$\dist_g : L_g\times L_g \to (0,+\infty)$ 
defined as follows
\begin{equation*}
\forall u,v \in L_g,
\ \dist_g(u,v) = \sum_{e\in\path_g(u,v)} \weight_g(e),
\end{equation*}
where $\path_g(u,v)$ is the set of edges on the path between
$u$ and $v$ in $T_g$.
We will refer to $\dist_g(u,v)$ as the {\em evolutionary
distance} between $u$ and $v$ under $g$.

For each pair of extant species $\{a,b\}$,
we compute the median 
$$
\distm(a,b) 
= \median
\{\dist_{g_i}(a,b) : i = 1,\ldots, N,\ \{a,b\} \subseteq L_{g_i}\}.
$$
We abort if a pair is 
not included in any of the gene trees.
We then use the distance matrix $\distm$ to 
build a tree using the Short Quartet Method~\cite{ErStSzWa:99a} (or any other 
statistically consistent, fast-converging distance-based
method). We will refer to this method as the
MedianTree (MT) method.
The algorithm is detailed in Figure~\ref{fig:plurality}.
\begin{figure*}[!ht]
\framebox{
\begin{minipage}{13cm}
{\small \textbf{Algorithm} MedianTree\\
\textit{Input:} $N$ alignments over the taxa $[n]$;\\
\textit{Output:} Estimated species phylogeny $\hat{T}$;

\begin{itemize}
\item For each gene $g_i$ and each pair of taxa $\{a,b\}$, 
compute the log-det distance $\widehat{\dist}_{g_i}(a,d)$.

\item For all pairs of taxa $\{a,b\}$, compute
$$
\widehat{\dist}_{\mathrm{m}}(a,b) 
= \median
\left\{\widehat{\dist}_{g_i}(a,b) : i = 1,\ldots, N,\ \{a,b\} \subseteq L_{g_i}\right\}.
$$

\item Using SQM~\cite{ErStSzWa:99a} on
the distance-matrix $\{\widehat{\dist}_{\mathrm{m}}(a,b)\}_{a,b \in [n]}$,
compute the tree $\hat{T}$ 
(or abort if no tree is found).

\item Output $\hat{T}$.

\end{itemize}

}
\end{minipage}
} \caption{
Algorithm MedianTree.
} \label{fig:plurality}
\end{figure*}

\paragraph{Probabilistic analysis}
Define the {\em maximum path weight (MPW)}
$$
\mrpw 
= \max\{\lgttotal_X : X \subseteq (L_s^+)^2\}.
$$
Then:
\begin{lemma}[Probability of a miss: Distance case]
\label{lem-missd}
Let $T_g = (V_g,E_g,L_g;\weight_g)$ be a gene tree distributed according 
to the random LGT model such that
$X = \{a,b\} \subseteq L_g$. 
Let $\dist_s(a,b)$ be the 
evolutionary distance between $a$ and $b$ under
the topology of the extant phylogeny
(that is, under the event that no LGT has occurred). Then
$$
\P[\dist_g(a,b) = \dist_s(a,b)  | X \subseteq L_g]
\geq \exp\left(-\mrpw\right).
$$
\end{lemma}
\begin{prevproof}{Lemma}{lem-missd}
The proof is similar to that of Lemma~\ref{lem-missq}.
\end{prevproof}

\begin{lemma}[Bound on path weight: Bounded-rates case]
\label{lem-relbrmd}
Under the Bounded-rates model, it holds that
$$
\mrpw = O\left(\maxlgt \log n^+\right).
$$
\end{lemma}
\begin{prevproof}{Lemma}{lem-relbrmd}
Note that
$$
\max\{\lgttotal_X : X \subseteq (L_s^+)^2\}
\leq 2 \maxlgt \maxtime \frac{\maxtime}{\mintime} \log_2 n^+.
$$
\end{prevproof}
\begin{lemma}[Bound on path weight: Yule case]
\label{lem-relyuled}
Under the Yule model, it holds that
$$
\mrpw = \Theta\left(\maxlgt \log n\right),
$$
with probability approaching $1$ as $n \to +\infty$.
\end{lemma}
\begin{prevproof}{Lemma}{lem-relyuled}
The proof is similar to that of Lemma~\ref{lem-relyule}.
\end{prevproof}

\begin{proof}(Theorems~\ref{thm:main1} 
and~\ref{thm:main2})
Using MT and Lemmas~\ref{lem-relbrmd} and~\ref{lem-relyuled},
the proof of Theorem~\ref{thm:main1} (and of Theorem~\ref{thm:main2}) follows
from the same lines as that of Theorem~\ref{thm:main1}.
Note however that our extra assumption on
the gene tree branch lengths is needed here 
to ensure that evolutionary distances are the same
across all genes.
\end{proof}

\subsection{Taking into account sequence length}

We have assumed so far that gene tree topologies
and evolutionary distances
are known perfectly. Of course, this is 
not the case in practice and the effect of
sequence length must be accounted for. 
One issue that arises is that LGT events may
create very short branches that are difficult
to infer. 
Nevertheless, we can prove the following.
We assume that sequence data is generated
independently on each gene tree according
to a GTR model. Evolutionary distances
are estimated using the log-det distance.
See e.g.~\cite{SempleSteel:03} for background
on GTR models of substitution and the log-det
distance.
We assume $n^- = 0$ for simplicity.
\begin{theorem}[Sequence-length requirements]
\label{thm:seq}
Under the Bounded-rates and Yule models for the species
phylogeny and the GTR model for
sequences, assuming that substitution rates
are bounded between constants,
a sequence length per gene polynomial in $n$
suffices
for the MT algorithm to succeed
if the number of genes is at most polynomial in $n$.
\end{theorem}
\begin{prevproof}{Theorem}{thm:seq}
We only discuss the Yule model. The argument
for the Bounded-rates model is similar.

In our second proof of Theorem~\ref{thm:main2}, 
we relied on the fact that, for every pair of taxa w.h.p.,
a strict majority of the gene tree evolutionary distances
{\em is not been affected by LGT.}
Hence, if the worst case estimation error on 
the evolutionary distances is $\eps$, then
the median of the estimated distances must be in
the interval $[\dist_s(a,b) - \eps, \dist_s(a,b) + \eps]$
for all pairs of taxa $a,b$.
Further, by the concentration bounds in~\cite{ErStSzWa:99b},
for the SQM step of our MT algorithm to return
the correct topology w.h.p., the sequence length
must scale as an exponential of the depth of the tree
divided by the square of the shortest branch length.

Under the Yule model, 
with probability approaching $1$, 
the depth of the tree is $O(\log n)$ 
(by the proof of Lemma~\ref{lem-relyule})
and
the shortest branch length (the minimum
of $O(n)$ exponentials with mean $O(1)$)
is $1/\poly(n)$.
Hence 
the result follows.\footnote{Note that unlike~\cite{ErStSzWa:99a} we use
the inter-speciation times generated by the
continuous-time branching process.
In particular their ``few logs'' result does not
apply to our setting.}
\end{prevproof}

\section{Discussion}

We have shown that a species phylogeny or network 
can be reconstructed despite high levels of random LGT 
and we have provided explicit quantitative bounds on tolerable rates of LGT. 
Moreover our analysis sheds light 
on effective approaches for species tree building in the presence of LGT. Several problems remain open: 
\begin{itemize}
\item
Galtier and Daubin~\cite{GaltierDaubin:08} 
hypothesize that random LGT only becomes a significant hurdle when the rate of LGT greatly exceeds the rate of diversification. In our setting this would imply that a value of $\lgttotal$ as high as $\Omega(n)$ may be achievable. 
Note that branches close to the leaves 
are particularly easy to reconstruct because
they lie on small quartet trees that are less likely
than deep ones to be hit by an LGT event. 
Is a recursive approach starting from
the leaves possible here? See~\cite{Mossel:04a,DaMoRo:11a} for recursive approaches in a related context.

\item
In a related problem, we have analyzed 
distance-based and 
quartet-based methods. 
A better understanding of bipartition-based
approaches is needed and may lead to a higher threshold for $\lgttotal$. 

\item What can be proved when a model
of extinction is incorporated?

\item What can be proved when the number of genes
is significantly less than $\log n$?

\item 
In the presence of highways, dealing with
more general network settings would be desirable.
Also our definition of highways as connecting
two edges is somewhat restrictive. In general,
one is also interested in preferential
genetic transfers between clades.

\item On the practical side, the predictions
made here should be further tested on real and
simulated datasets. 
We note that there is extisting work
in this direction~\cite{BeHaRa:05,GeWaKi:05,Galtier:07,PuWoKo:09,PuWoKo:10,KoPuWo:11,BaBaGo+:11}.
\end{itemize}

\bibliographystyle{alpha}
\bibliography{own,thesis,RECOMB12}

\clearpage

\appendix

\end{document}